\documentclass[%
final
]{dmtcs-episciences}

\pdftrailerid{}
\usepackage[utf8]{inputenc}
\usepackage{subfigure}

\usepackage[round]{natbib}

\def\({\left(}
\def\){\right)}  
\def\ee{{\rm e}}
\def\llceil{\left\lceil}
\def\rrceil{\right\rceil}  

\let\epsilon\varepsilon

\newcommand{\cA}{\mathcal{A}}
\newcommand{\cH}{\mathcal{H}}
\newcommand{\cI}{\mathcal{I}}
\newcommand{\cL}{\mathcal{L}}
\newcommand{\cP}{\mathcal{P}}
\newcommand{\cR}{\mathcal{R}}
\newcommand{\cU}{\mathcal{U}}
\newcommand{\EE}{\mathbb{E}}
\newcommand{\PP}{\mathbb{P}}
\newcommand{\RR}{\mathbb{R}}
\newcommand{\opt}{{\rm OPT}}
\def\HAR{\textsc{Harmonic}}
\def\val{\mathop{{\rm w}}\nolimits}

\usepackage{amsthm}
\newtheorem{theorem}              {Theorem}
\newtheorem{lemma}      [theorem] {Lemma}

\newtheorem{fact}       [theorem] {Fact}
\theoremstyle{definition}
\newtheorem{remark}     [theorem] {Remark}
\newtheorem{definition} [theorem] {Definition}

\usepackage{tikz}
\usepackage{amsmath}
\usepackage{enumerate}

\author{%
  Yoshiharu Kohayakawa\affiliationmark{1}\thanks{Partially supported
    by CNPq (311412/2018-1, 423833/2018-9) and FAPESP (2018/04876-1).}
  \and
  Fl\'avio Keidi Miyazawa\affiliationmark{2}\thanks{Partially supported
    by CNPq (314366/2018-0, 425340/2016-3) and FAPESP (2015/11937-9,
  2016/01860-1).}\\ 
  \and
  Yoshiko Wakabayashi\affiliationmark{1}\thanks{Partially supported by CNPq
  (306464/2016-0, 423833/2018-9) and FAPESP (2015/11937-9).  
  FAPESP is the S\~ao Paulo Research Foundation.  CNPq is the National
  Council for Scientific and Technological Development of Brazil.}}
\title
[Online bounded space hypercube bin packing]{A tight lower bound for
  the online bounded space hypercube bin packing problem\thanks{This
  research was partially supported by CAPES (Finance Code 001).}}  
\affiliation{%
  Institute of Mathematics and Statistics, University of S\~ao Paulo, Brazil\\
  Institute of Computing, University of Campinas, Brazil
}
\received{2021-07-30}
\accepted{2021-08-20}
\begin{document}
\publicationdetails{23}{2021}{3}{7}{8325}
\maketitle
\begin{abstract}
  In the $d$-dimensional hypercube bin packing problem, a given
  list of $d$-dimensional hypercubes must be packed into the
  smallest number of hypercube bins. Epstein and van Stee~[SIAM
  J. Comput.~35 (2005)] showed that the asymptotic performance
  ratio~$\rho$ of the online bounded space variant is $\Omega(\log d)$
  and $O(d/\log d)$, and conjectured that it is $\Theta(\log d)$. We
  show that~$\rho$ is in fact $\Theta(d/\log d)$, using probabilistic
  arguments.
\end{abstract}

\section{Introduction}
\label{sec:introduction}

The bin packing problem is an iconic problem in combinatorial
optimization, which has been investigated intensively from many
different viewpoints.  In particular, it has served as a proving
ground for new approaches to the development and analysis of
approximation and online algorithms, as well as for the development of
average case analysis techniques~(see \cite{CoffmanGJ97,CoffmanCGMV13}).

We prove a lower bound for a variant of the bin packing problem, in
which the items to be packed are $d$-dimensional hypercubes, also
referred to as \emph{$d$-hypercubes} or simply hypercubes, when the
dimension is clear.  More precisely, we prove a tight lower bound for
the \emph{online bounded space $d$-hypercube bin packing problem},
settling an open problem raised by~\cite{EpsteinS05}.  
Before we state our result (Theorem~\ref{thm:lwbd_prbsa}), we
introduce the required concepts and definitions and discuss briefly
the relevant literature.

The \emph{$d$-hypercube bin packing problem} ($d$-CPP) is defined as
follows.  We are given a list $L$ of \textit{items}, where each
item~$h$ in~$L$ is a $d$-hypercube of side length $s(h)\leq1$, and an
unlimited number of \textit{bins}, each of which is a unit
$d$-hypercube (that is, a $d$-hypercube of side length~$1$).
The goal is to find a packing~$\cP$ of the items in~$L$ into the
smallest possible number of bins. 
More precisely, we have to assign each item $h$ to a bin, and specify
its position
in that bin.
We require that the items be placed parallel to the axes of the bin
and, crucially, we require that the items in a bin should not overlap.
The \emph{size}~$|\cP|$ of the packing~$\cP$ is the number of
\emph{used} bins (those with assigned items).

The $d$-hypercube bin packing problem ($d$-CPP) is in fact a special
case of the \textit{$d$-dimensional bin packing problem} ($d$-BPP), in
which one has to pack $d$-dimensional parallelepipeds into
$d$-dimensional unit bins.  For
$d=1$, both problems reduce to the well known \textit{bin packing
  problem}. 

In the \emph{online} variant of $d$-CPP, the items arrive sequentially
and each item must be placed in some bin as soon as it arrives,
without knowledge of the next items.  The \emph{online bounded space}
variant of $d$-CPP is a restricted variant of online $d$-CPP.
Whenever a new empty bin is used in the packing process, it is
considered an \textit{open} bin and it remains so until it is
considered \textit{closed}, after which point it is not allowed to
accept other items.  In this variant, regardless of the
instance~$I$, at every point of the process, not more than~$M$ bins
should be open, where~$M$ is some constant that does not depend
on~$I$.

As usual for bin packing problems, we consider the asymptotic
performance ratio to measure the quality of algorithms.  For an
algorithm~$\cA$ and an input list~$L$, let~$\cA(L)$ be the number of bins
used by the solution produced by~$\cA$ for the list~$L$.  Furthermore,
let~$\opt(L)=\min|\cP|$, where the minimum is taken over all
possible packings~$\cP$ of~$L$ into unit bins.
The \emph{asymptotic performance ratio} of~$\cA$ is
\begin{equation}
\cR_{\cA}^{\infty} = \limsup_{n\rightarrow\infty}
               \sup_L\left\{  \frac{\cA(L)}{\opt(L)}: \opt(L)=n\right\}.
\end{equation}
Given a packing problem $\Pi$, the \emph{optimal asymptotic performance
  ratio} for $\Pi$ is
\begin{equation}
  \cR_{\Pi}^{\infty}=\inf\left\{\cR_{\cA}^{\infty}: \cA\mbox{ is an
    algorithm for }\Pi\right\}.
\end{equation}

Many results have been obtained for online $d$-BPP and $d$-CPP
(see, \textit{e.g.}, \cite{%
BaloghBDEL19,  %
BaloghBG12,         %
HeydrichS2016a,        %
christensen17:_approx, %
HeydrichS16b,          %
Seiden02,
Vliet92%
}).  
In our brief discussion of the literature below, we restrict ourselves
to the online bounded space versions of $d$-BPP and $d$-CPP.

For online bounded space $1$-BPP, \cite{lee85:_simpl_online} gave an
algorithm called~$\HAR_M$ with asymptotic performance ratio at
most~$(1+\epsilon)\Pi_\infty$, where~$\epsilon\to0$ as~$M\to\infty$,
and~$\Pi_\infty\approx1.69103$ is a certain explicitly defined
constant.  These authors also proved that no algorithm for online
bounded space $1$-BPP can have asymptotic performance ratio smaller
than~$\Pi_\infty$.  For online bounded space $d$-BPP for general~$d$,
a lower bound of~$\Pi_{\infty}^d$ was implicitly proved by
\cite{CsirikV93}, and \cite{EpsteinS05} proved an asymptotically
matching upper bound.

For online bounded space $d$-CPP, \cite{EpsteinS05} proved that its
asymptotic performance ratio is~$\Omega(\log d)$ and~$O(d/\log d)$,
and conjectured that it is $\Theta(\log d)$.  They also gave an
optimal algorithm for this problem, but left as an interesting open
problem to determine its asymptotic performance ratio.
Later, \cite{EpsteinS07} gave lower and upper bounds for
$d\in\{2,\ldots,7\}$.

Our main contribution is an $\Omega(d/\log d)$ lower bound for online
bounded space $d$-CPP.  In view of previous results by
\cite{EpsteinS05}, we obtain that the asymptotic performance ratio of
this problem is~$\Theta(d/\log d)$, settling an open problem posed by
those authors.  To prove our lower bound, we follow a well known
approach (see ~\cite{lee85:_simpl_online} and \cite{yao80:_new}), which
requires the proof of the existence of a packing with a suitably large
`weight', for a certain definition of weight.  The novelty here is
that we prove the existence of such a packing with the probabilistic
method.

To conclude this section, we mention that the technique that we
present here may also be used to obtain lower bounds for the prices of
anarchy of a game theoretic version of $d$-CPP, called \emph{selfish
  $d$-hypercube bin packing game}.  As this topic requires the
introduction of a number of concepts, we just mention the main results
for readers familiar with this line of research: \textsl{for every
  large enough~$d$, the asymptotic price of anarchy
  {\rm(}respectively, strong price of anarchy{\rm)} of the selfish
  $d$-hypercube bin packing game is $\Omega(d/\log d)$
  {\rm(}respectively, $\Omega(\log d)${\rm)}}.  The proof of one of
the results can be found in~\cite{KohayakawaMW2017}.  A preliminary
version of this work (\cite{KohayakawaMW18}) appeared in the
proceedings of LATIN~2018.

\section{Notation and homogeneous packings}
\label{sec:preliminaries}

The $d$-hypercubes~$Q_k^+=Q_k^d(\epsilon)$ defined below will be
important in what follows.

\begin{definition}
  \label{def:Q_k^d(gepsilon)}
  Let~$d\geq2$ be an integer.  For every integer~$k\geq2$
  and~$\epsilon>0$, let
  \begin{equation}
    \label{eq:Q_k^d(epsilon).def}
    Q_k^+=Q_k^d(\epsilon)=\(0,{1+\epsilon\over k}\)^d
    =\left\{x\in\RR\colon0<x<{1+\epsilon\over k}\right\}^d\subset\RR^d
  \end{equation}
  be the open $d$-hypercube of side length~$(1+\epsilon)/k$
  `based' at the origin.
\end{definition}

\subsection{Homogeneous packings}
\label{sec:homogeneous}
We shall be interested in certain types of packings %
of hypercubes into a unit bin.

\begin{definition}[Homogeneous packings~$\cH_k^+=\cH_k^d(\epsilon)$]
  \label{def:cH}
  Let~$d\geq2$ be fixed.  For any integer $k\geq2$
  and~$0<\epsilon\leq1/(k-1)$, a packing of $(k-1)^d$ copies
  of~$Q_k^+=Q_k^d(\epsilon)$ into a unit bin is said to be a
  \textit{packing of type~$\cH_k^+=\cH_k^d(\epsilon)$}.  Packings of
  type~$\cH_k^+$ will be called \textit{homogeneous packings}.
\end{definition}

In the definition above, the upper bound on~$\epsilon$ guarantees
that~$(k-1)^d$ copies of~$Q_k^+$ \textit{can} be packed into
a unit bin (and hence~$\cH_k^+$ exists): it suffices to
note that, under that assumption on~$\epsilon$, we have
$(k-1)(1+\epsilon)/k\leq1$.
Homogeneous packings are important because
they can be used to create instances for which any bounded space
algorithm performs badly (see \cite{EpsteinS05,EpsteinS07}). 

\section{The central lemma and the main theorem}
\label{sec:packing}
The key result used in the proof of our main theorem
(Theorem~\ref{thm:lwbd_prbsa}) is 
the existence of a certain packing, stated in
Lemma~\ref{lem:large_value_pack} below.  Since this lemma is somewhat
technical, we first informally describe a related result.

Consider the~$S-1$ homogeneous packings~$\cH_k^+$ ($k=2,\dots,S$),
where~$S=\lceil cd/\log d\rceil $ for a small positive constant~$c$.
Suppose also that~$0<\epsilon\leq\epsilon_0(d)$ for some
small~$\epsilon_0(d)$.  Suppose we assemble a list~$\cI$ of
$d$-hypercubes from these~$S-1$ homogeneous packings~$\cH_k^+$ by
selecting~$90\%$ of the members of each such~$\cH_k^+$.  The following
holds: (*)~\textsl{there is a packing of~$\cI$ into a \textbf{single}
  unit bin as long as~$d$ is sufficiently large}.  This fact is behind
the proof of our central lemma, Lemma~\ref{lem:large_value_pack},
stated in what follows. Fact~(*) might look surprising at first sight, as the
homogeneous packings~$\cH_k^+$ appear to have reasonably high
occupancy.

We now give some definitions needed for the statement of
Lemma~\ref{lem:large_value_pack}.

\begin{definition}[$\epsilon$-packings]
  \label{def:(1+epsilon)/ZZ}
  A packing~$\cU$ of $d$-hypercubes into a unit bin is called an
  \textit{$\epsilon$-packing} if, for every member~$Q$ of~$\cU$, there
  is some integer~$k\geq2$ such that~$Q$ is a copy
  of~$Q_k^+=Q_k^d(\epsilon)$.
\end{definition}

Let~$\cU$ be an $\epsilon$-packing for some~$\epsilon>0$.  Let
\begin{equation}
  \label{eq:K(U)_def}
  K(\cU)=\{k\geq2\colon\cU\text{ contains a copy of }Q_k^+\}.
\end{equation}
For every~$k\in K(\cU)$, let 
\begin{equation}
  \label{eq:nu_k(U)_def}
  \nu_k(\cU)\text{ be the total number of copies of
  }Q_k^+\text{ in }\cU.
\end{equation}
Clearly, we have~$0\leq\nu_k(\cU)\leq(k-1)^d$ for every~$k$ (recall
that~$\epsilon>0$).  Finally, we define the \textit{weight} of~$\cU$
as
\begin{equation}
  \label{eq:val_def}
  \val(\cU)=\sum_{k\in K(\cU)}{\nu_k(\cU)\over(k-1)^d}.
\end{equation}

We shall be interested in $\epsilon$-packings~$\cU$ with large weight.
Our main lemma is as follows.

\begin{lemma}[Central lemma] \label{lem:large_value_pack}
  There is an absolute constant~$d_0$ for which the following holds
  for any~$d\geq d_0$.  For any~$0<\epsilon\leq d^{-2}$, the unit bin
  admits an $\epsilon$-packing~$\cU$ with
  \begin{equation}
    \label{eq:large_value_pack}
    \val(\cU)\geq{d\over5\log d}.
  \end{equation}
\end{lemma}

In~\eqref{eq:large_value_pack} and in what follows, $\log x$ stands
for the natural logarithm of~$x$.  The proof of
Lemma~\ref{lem:large_value_pack} is postponed to
Section~\ref{sec:packing_proof}.  We now deduce our main result, 
Theorem~\ref{thm:lwbd_prbsa}, from Lemma~\ref{lem:large_value_pack},
following the approach used by \cite{lee85:_simpl_online}.  For experts in the
area, given Lemma~\ref{lem:large_value_pack}, the proof of
Theorem~\ref{thm:lwbd_prbsa} is routine.  The short proof below is
included for the benefit of non-experts.

\begin{theorem}[Main Theorem]
  \label{thm:lwbd_prbsa}
  There is an absolute constant~$d_0$ such that, for any~$d\geq d_0$,
  the asymptotic performance ratio of the online bounded space
  $d$-hypercube bin packing problem is at least~$d/10\log d$.
\end{theorem}

\begin{proof} %
  Let~$\cA$ be any algorithm for the online bounded space
  $d$-hypercube bin packing problem.  Let~$M$ be the maximum number of
  bins that~$\cA$ leaves open during its execution.  To
  prove that~$\cA$ has asymptotic performance ratio at
  least~$d/10\log d$ if~$d$ is large enough, we construct a suitable
  instance~$\cI$ for~$\cA$.

  Let~$d_0$ be as in Lemma~\ref{lem:large_value_pack} and
  suppose~$d\geq d_0$.  Fix any~$\epsilon$
  with~$0<\epsilon\leq d^{-2}$ and let~$\cU$ be a packing as in the
  statement of Lemma~\ref{lem:large_value_pack}.  The instance~$\cI$
  will be constructed as follows.  First, we choose a suitable
  integer~$N$ and take~$2MN$ copies of~$\cU$.  We then construct~$\cI$
  by arranging the hypercubes in these copies in a linear order, with
  all the hypercubes of the same size appearing together.  Let us now
  formally describe~$\cI$.

Let~$N=\prod_{k\in K(\cU)}(k-1)^d$.
Recall that~$\cU$ contains~$\nu_k(\cU)$ copies of~$Q_k^+$
for every~$k\in K(\cU)$.  Let~$K=|K(\cU)|$ and
suppose~$K(\cU)=\{k_1,\dots,k_K\}$.  The instance~$\cI$ that we shall
construct is the concatenation of~$K$ segments,
say~$\cI=\cI_1\dots\cI_K$,
with each segment~$\cI_\ell$ ($1\leq\ell\leq K$) composed of a
sequence of~$f(\ell)=2MN\nu_{k_\ell}(\cU)$ copies
of~$Q_{k_\ell}^+$. 
This completes the definition of our instance~$\cI$.

The following assertion, to be used later, concerning the offline
packing of the hypercubes in~$\cI$ is clear, as we obtained~$\cI$ by
rearranging the hypercubes in~$2MN$ copies of~$\cU$.
\begin{equation}
  \label{eq:1}
  \text{The hypercubes in~$\cI$ can be packed into at
    most~$2MN$  unit bins.} 
\end{equation}

We now prove that, when~$\cA$ is given the instance~$\cI$ above, it
will have performance ratio at least as bad as~$\val(\cU)/2$.  In view
of~\eqref{eq:large_value_pack} in Lemma~\ref{lem:large_value_pack},
this will complete the proof of Theorem~\ref{thm:lwbd_prbsa}.

Let us examine the behaviour of~$\cA$ when it is given input~$\cI$.
Fix~$1\leq\ell\leq K$ and suppose that~$\cA$ has already seen the 
hypercubes in~$\cI_1\dots\cI_{\ell-1}$ and it has already packed them
somehow.  We now consider what happens when~$\cA$ examines
the~$f(\ell)$ hypercubes in~$\cI_\ell$, which are all copies
of~$Q_{k_\ell}^+$.

Clearly, since~$\epsilon>0$, the~$f(\ell)$ copies
of~$Q_{k_\ell}^+$ in~$\cI_\ell$ cannot be packed into fewer
than
\begin{equation}
  \label{eq:cI_ell_lwbd}
  {f(\ell)\over(k_\ell-1)^d}
  ={2MN\nu_{k_\ell}(\cU)\over(k_\ell-1)^d}
  \geq{MN\nu_{k_\ell}(\cU)\over(k_\ell-1)^d}+M
\end{equation}
unit bins.  Therefore, even if some hypercubes in~$\cI_\ell$ are placed in
bins still left open after the processing of~$\cI_1\dots\cI_{\ell-1}$,
the hypercubes in~$\cI_\ell$ will add at
least~$MN\nu_{k_\ell}(\cU)/(k_\ell-1)^d$ new bins to the current
output of~$\cA$. 
Thus, the total number of bins that~$\cA$ will use when
processing~$\cI$ is at least
\begin{equation}
  \label{eq:cA_lwbd}
  \sum_{k\in K(\cU)}{MN\nu_k(\cU)\over(k-1)^d}
  =MN\sum_{k\in K(\cU)}{\nu_k(\cU)\over(k-1)^d}
  =MN\val(\cU).
\end{equation}
In view of~\eqref{eq:1}, it follows that the asymptotic
performance ratio of~$\cA$ is at least
\begin{equation}
  \label{eq:cA_apr_lwbd}
  {MN\val(\cU)\over2MN}={1\over2}\val(\cU),
\end{equation}
as claimed.  This completes the proof of
Theorem~\ref{thm:lwbd_prbsa}. 
\end{proof}

\section{Proof of Lemma~\ref{lem:large_value_pack}}
\label{sec:packing_proof}

The $\epsilon$-packing~$\,\cU$ whose existence is asserted in our
central lemma, Lemma~\ref{lem:large_value_pack}, will be described in
terms of certain `codes', that is, sets of `codewords' or simply
`words'.  We shall use such codes to `place' copies of
certain~$Q_k^+=Q_k^d(\epsilon)$ in the packing~$\cU$.  We make this
precise in Section~\ref{sec:plac-hyperc-accord}.  The proof of
the existence of appropriate codes will be given in
Section~\ref{sec:words_fam}.  The proof of
Lemma~\ref{lem:large_value_pack} is given in
Section~\ref{sec:packing-cu_epsilon}. 

\subsection{Placing hypercubes according to codewords}
\label{sec:plac-hyperc-accord}
Let~$d$ and~$k\geq2$ be fixed.  Let a \textit{$d$-letter
word}~$w\in[k]^d$ from the \textit{alphabet}~$[k]=\{1,\dots,k\}$ be
given.  In what follows, we shall fix
some~$0<\epsilon\leq\epsilon_0(d)$ and we shall consider
translations~$Q(w)=Q^{(k)}(w)$ of the hypercube~$Q_k^+$
specified by such words~$w$ in a certain way (for the definition
of~$Q_k^+=Q_k^d(\epsilon)$, recall~\eqref{eq:Q_k^d(epsilon).def}).
Furthermore, later, we shall consider certain
\textit{sets}~$L_k\subset[k]^d$ 
of such words and we shall define packings of the
form~$\cP_{L_k}=\{Q(w)\colon w\in L_k\}$.  Note that~$\cP_{L_k}$ is
composed of copies of~$Q_k^+$.  To obtain the packing~$\cU$
whose existence is asserted in Lemma~\ref{lem:large_value_pack}, we
shall consider the union of such packings~$\cP_{L_k}$ for
$k=2,\dots,S$, with~$S=\lceil2d/9\log d\rceil$
and certain families~$\cL=\{L_k\colon2\leq k\leq S\}$
(see Lemma~\ref{lem:sep_langs}).

\begin{figure}[h]
  \centering
  \begin{tikzpicture}
    \draw[line width=1.5pt,-stealth] (-.3,0) -- (6.5,0)
    node[anchor=west] {$i$};
    \draw[line width=1.5pt,-stealth] (0,-.3) -- (0,6.5)
    node[anchor=south] {$j$};
    
    \foreach \i in {1,...,5} {
      \path (\i*1.05,0) coordinate (X\i);
      \path (\i*1.05,6) coordinate (X\i\i);
      \path (0,\i*1.05) coordinate (Y\i);
      \path (6,\i*1.05) coordinate (Y\i\i);
      \draw[dotted] (X\i) -- (X\i\i);
      \draw[dotted] (Y\i) -- (Y\i\i);
    }
    \draw (0,0) node[anchor=east] {\footnotesize$x^{(k)}(1)=0$\quad\null};
    \draw (Y1) node[anchor=east] {\footnotesize$x^{(k)}(2)=y^{(k)}(1)$};
    \draw (Y11) node[anchor=west] {\footnotesize$(1+\epsilon)/k$};
    \draw (Y2) node[anchor=east] {\footnotesize$x^{(k)}(3)=y^{(k)}(2)$};
    \draw (Y22) node[anchor=west] {\footnotesize$2(1+\epsilon)/k$};
    \draw (Y4) node[anchor=east] {\footnotesize$x^{(k)}(k-1)=y^{(k)}(k-2)$};
    \draw (Y44) node[anchor=west] {\footnotesize$(k-2)(1+\epsilon)/k$};
    \draw[->] (Y5) +(-1.95,0) -- (-.1,5*1.05);
    \draw (Y5) +(-1.95,0) node[anchor=east] {\footnotesize$y^{(k)}(k-1)$};
    \draw[->] (8.3,5*1.05) -- (6.1,5*1.05);
    \draw (8.3,5*1.05) node[anchor=west] {\footnotesize$(k-1)(1+\epsilon)/k$};
    \path (6,0) coordinate (X6);
    \path (6,6) coordinate (X66);
    \path (0,6) coordinate (Y6);
    \path (6,6) coordinate (Y66);
    \draw (Y6) node[anchor=east] {\footnotesize$y^{(k)}(k)=1$};
    \draw (Y66) node[anchor=west] {\footnotesize$1$};
    \draw[dotted] (X6) -- (X66);
    \draw[dotted] (Y6) -- (Y66);
    \path (4.95,0) coordinate (X5-);
    \path (4.95,6) coordinate (X55-);
    \path (0,4.95) coordinate (Y5-);
    \path (6,4.95) coordinate (Y55-);
    \draw (Y5-) node[anchor=east] {\footnotesize$x^{(k)}(k)$};
    \draw (Y55-) node[anchor=west] {\footnotesize$1-(1+\epsilon)/k$};
    \draw[dotted] (X5-) -- (X55-);
    \draw[dotted] (Y5-) -- (Y55-);

    \path (1.05,1.05) coordinate (A);
    \path (2.1,1.05) coordinate (B);
    \path (2.1,2.1) coordinate (C);
    \path (1.05,2.1) coordinate (D);
    \draw[fill,color=black!30,draw=black,
    line width=2pt,opacity=0.5]
    (A) -- (B) -- (C) -- (D) -- cycle;
    \draw (1.5*1.05,1.5*1.05) node{\footnotesize$Q(w)$};
    
    \path (0,4.95) coordinate (AA);
    \path (1.05,4.95) coordinate (BB);
    \path (1.05,6) coordinate (CC);
    \path (0,6) coordinate (DD);
    \draw[fill,color=black!30,draw=black,
    line width=2pt,opacity=0.5]
    (AA) -- (BB) -- (CC) -- (DD) -- cycle;
    \draw (.5*1.05,5.48) node{\footnotesize$Q(w')$};
    
    \path (2*1.05,4*1.05) coordinate (AAA);
    \path (3*1.05,4*1.05) coordinate (BBB);
    \path (3*1.05,5*1.05) coordinate (CCC);
    \path (2*1.05,5*1.05) coordinate (DDD);
    \draw[fill,color=black!30,draw=black,
    line width=2pt,opacity=0.5]
    (AAA) -- (BBB) -- (CCC) -- (DDD) -- cycle;
    \draw (2.5*1.05,4.5*1.05) node{\footnotesize$Q(w'')$};

    \path (4*1.05,5*1.05) coordinate (A4);
    \path (5*1.05,5*1.05) coordinate (B4);
    \path (5*1.05,6*1.05) coordinate (C4);
    \path (4*1.05,6*1.05) coordinate (D4);
    \draw[fill,color=black!30,draw=black,
    line width=2pt,opacity=0.5]
    (A4) -- (B4) -- (C4) -- (D4) -- cycle;
    \draw (4.5*1.05,5.5*1.05) node{\footnotesize$Q'$};

  \end{tikzpicture}
  \caption{Projections on the $(i,j)$-plane of hypercubes~$Q(w)$, 
    $Q(w')$ and~$Q(w'')$ with $w_i=w_j=2$, $w_i'=1$ and~$w_j'=k$,
    and~$w''_i=3$ and~$w_j''=k-1$.  The hypercube~$Q'$ is not contained
    in~$[0,1]^d$.} 
  \label{fig:Qw}
\end{figure}
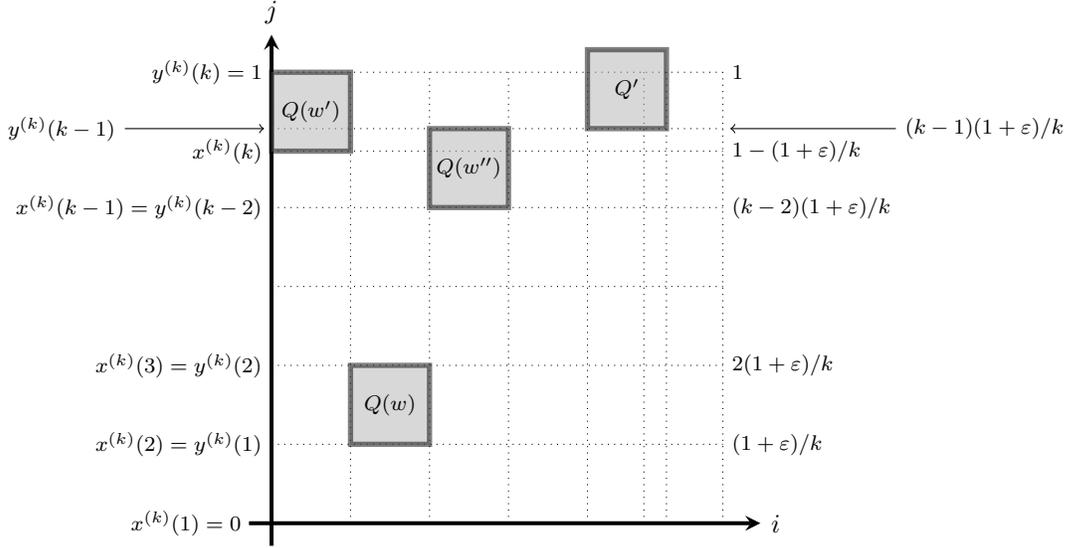

Let us now define~$Q(w)=Q^{(k)}(w)$, the translation of~$Q_k^+$
specified by~$w=(w_1,\dots,w_d)\in[k]^d$.
For~$w=(w_i)_{1\leq i\leq d}$ with~$w_i<k$ for every~$i$, we
let~$Q(w)$ be the translation
\begin{equation}
  \label{eq:2}
  x[w]+Q_k^+=\{x[w]+z\colon z\in Q_k^+\}
\end{equation}
of~$Q_k^+$, where
\begin{equation}
  \label{eq:3}
  x[w]={1+\epsilon\over k}(w_1-1,\dots,w_d-1).
\end{equation}
Thus, while~$Q_k^+$ has its `base point' at the origin, $Q(w)$~has its
base point at~$x[w]$ (see~$Q(w)$ and~$Q(w'')$ in Figure~\ref{fig:Qw}).

In what follows, we shall always have~$0<\epsilon<1/(k-1)$.
Therefore, if~$w_i<k$ for every~$i$, then~$Q(w)$ is contained in the
unit hypercube~$[0,1]^d$, whereas if~$w_i=k$ for some~$i$, then
$x[w]+Q_k^+=\{x[w]+z\colon z\in Q_k^+\}$ with~$x[w]$ as defined
in~\eqref{eq:3} is \textit{not} contained in~$[0,1]^d$ (see~$Q'$ in
Figure~\ref{fig:Qw}).  Since we want~$Q(w)$ to be contained
in~$[0,1]^d$ for every~$w\in[k]^d$, we actually define~$x[w]$ as
in~\eqref{eq:7} below.

\begin{definition}[Base point coordinates of~$Q(w)$]
  \label{def:x^k(j)}
  For every~$k\geq2$ and~$0<\epsilon<1/(k-1)$, let
  \begin{equation}
    \label{eq:x^q(j)_def}
    x^{(k)}(v)=x_\epsilon^{(k)}(v)=
    \begin{cases}\displaystyle
      {(1+\epsilon)(v-1)\over k},
      &\text{if }1\leq v<k,\\[2.5ex]
      \displaystyle
      1-{1+\epsilon\over k},
      &\text{if }v=k.
    \end{cases}
  \end{equation}
  For~$w=(w_1,\dots,w_d)\in[k]^d$, let
  \begin{equation}
    \label{eq:7}
    x[w]=(x^{(k)}(w_1),\dots,x^{(k)}(w_d)).
  \end{equation}
  Finally, for convenience, for~$1\leq v\leq k$, let
  \begin{equation}
    \label{eq:y^{(q)}(j)}
    y^{(k)}(v)=x^{(k)}(v)+{1+\epsilon\over k}.
  \end{equation}
\end{definition}

We now state three simple facts that the reader may find useful to
check on their own to get used to the definitions above.  First, note
that $\cP=\{Q(w)\colon w\in[k-1]^d\}$ is a packing of~$(k-1)^d$ copies
of~$Q_k^+$ into the unit bin~$[0,1]^d$; that is, $\cP$~is a packing of
type~$\cH_k^+$ (recall Definition~\ref{def:cH}).  Secondly,
$\{Q(w)\colon w\in[k]^d\}$ is \textit{not} a packing.  Finally,
$\{Q(w)\colon w\in[k]^d\text{ with }w_i\neq k-1\text{ for every }i\}$
\textit{is} a packing (and is also a packing of type~$\cH_k^+$).

Note that, because~$\epsilon<1/(k-1)$, for every~$k\geq2$, we have
\begin{multline}
  \label{eq:x^d(j).2}
  \qquad
  0=x^{(k)}(1)<y^{(k)}(1)=x^{(k)}(2)<y^{(k)}(2)=x^{(k)}(3)<\cdots
  <y^{(k)}(k-2) \\ 
  =x^{(k)}(k-1)
  <x^{(k)}(k)<y^{(k)}(k-1)<y^{(k)}(k)=1\qquad
\end{multline}
(see Figure~\ref{fig:Qw}).  For every~$k\geq2$ and~$1\leq v\leq k$,
let
\begin{equation}
  \label{eq:I^(k)(j)_def}
  I^{(k)}(v)=(x^{(k)}(v),y^{(k)}(v))\subset[0,1].
\end{equation}
Finally, note that 
\begin{equation}
  \label{eq:Q(w)_also}
  Q(w)=Q^{(k)}(w)=x[w]+Q_k^+
  =I^{(k)}(w_1)\times\dots\times I^{(k)}(w_d)
  \subset[0,1]^d.
\end{equation}
We close this section observing the following.

\begin{fact}    %
  \label{fact:gap}
  The following assertions hold for any positive~$S$. 
  \begin{enumerate}[(i)]
  \item\label{enum:fact_gap.i} Suppose~$2\leq k<k'\leq S$
    and~$0<\epsilon\leq S^{-2}$.  Then
    \begin{equation}
      \label{eq:gap}
      y^{(k)}(k-1)<x^{(k')}(k').
    \end{equation}
    In particular, the intervals~$I^{(k)}(v)$ $(1\leq v<k)$ are disjoint
    from~$I^{(k')}(k')$.
  \item\label{enum:fact_gap.ii} For any~$2\leq k\leq S$, the
    intervals~$I^{(k)}(v)$ $(1\leq v\leq k)$ are pairwise disjoint,
    except for the single pair formed by~$I^{(k)}(k-1)$
    and~$I^{(k)}(k)$.
  \end{enumerate}
\end{fact}
\begin{proof}
  Assertion~\eqref{enum:fact_gap.ii} is clear
  (recall~\eqref{eq:x^d(j).2}).  The second assertion
  in~\eqref{enum:fact_gap.i} follows from inequality~\eqref{eq:gap}, and
  therefore it suffices to verify that inequality.  We
  have~$y^{(k)}(k-1)=x^{(k)}(k-1)+(1+\epsilon)/k=(k-1)(1+\epsilon)/k
  =1+\epsilon-(1+\epsilon)/k$.  Moreover,
  $x^{(k')}(k')=1-(1+\epsilon)/k'$.  Therefore,
  \eqref{eq:gap}~is equivalent to
  \begin{equation}
    \label{eq:gap.1}
    \epsilon<(1+\epsilon)\({1\over k}-{1\over k'}\).
  \end{equation}
  Since~$k+1\leq k'\leq S$ and~$\epsilon\leq S^{-2}$,
  inequality~\eqref{eq:gap.1} does hold.
\end{proof}

\subsection{Separated families of gapped codes}
\label{sec:words_fam}
Let an integer~$d\geq2$ be fixed.  We shall consider sets of
words~$L_k\subset[k]^d=\{1,\dots,k\}^d$ for~$k\geq2$.  We refer to
such~$L_k$ as \textit{codes} or $k$-\textit{codes}.  As discussed in
the beginning of Section~\ref{sec:plac-hyperc-accord}, we shall design
such~$L_k$ to specify packings $\cP_{L_k}=\{Q(w)\colon w\in L_k\}$.

We start with the following definition.

\begin{definition}[Gapped codes]
  \label{def:gapped_langs}
  Suppose~$k\geq2$ and let a $k$-code $L_k\subset[k]^d$ be given.
  We say that~$L_k$ \textit{misses}~$j$ at coordinate~$i_0$ if every
  word~$w=(w_i)_{1\leq i\leq d}$ in~$L_k$ is such that~$w_{i_0}\neq
  j$.  Furthermore, $L_k$~is said to be \textit{gapped} if, for
  each~$1\leq i\leq d$, either~$L_k$ misses~$k-1$ at~$i$ or~$L_k$
  misses~$k$ at~$i$.
\end{definition}

Suppose~$L_k$ is a gapped code, and suppose~$w=(w_i)_{1\leq i\leq d}$
and~$w'=(w_i')_{1\leq i\leq d}$ are distinct words in~$L_k$.
Then~$Q(w)$ and~$Q(w')$ do not overlap: this can be checked
from~\eqref{eq:Q(w)_also} and
Fact~\ref{fact:gap}\eqref{enum:fact_gap.ii}.
Thus, if~$L_k$ is gapped, then
\begin{equation}
  \label{eq:cP_L_k_def}
  \cP_{L_k}=\{Q(w)\colon w\in L_k\}
\end{equation}
is a packing.

We now introduce a certain notion of `compatibility' between two
codes~$L_k$ and~$L_{k'}$, so that~$\cP_{L_k}$ and~$\cP_{L_{k'}}$ can be
put together to obtain a packing if they come from `compatible'
codes~$L_k$ and~$L_{k'}$.

\begin{definition}[Separated codes]
  \label{def:separated_langs}
  Suppose~$2\leq k<k'$ and~$L_k\subset[k]^d$ and~$L_{k'}\subset[k']^d$
  are given.  We say that~$L_k$ and~$L_{k'}$ are \textit{separated}
  if, for any~$w=(w_i)_{1\leq i\leq d}\in L_k$ and
  any~$w'=(w_i')_{1\leq i\leq d}\in L_{k'}$, there is some~$i$ such
  that~$w_i<k<k'=w_i'$.  
\end{definition}

Suppose~$L_k$ and~$L_{k'}$ are gapped and separated and
suppose~$k<k'\leq S$ and~$\epsilon\leq S^{-2}$ for some~$S$ (we shall
later set~$S$ to be a certain value~$S(d)$).  Consider the
packings~$\cP_{L_k}$ and~$\cP_{L_{k'}}$ as defined in~\eqref{eq:cP_L_k_def}.
Fact~\ref{fact:gap}\eqref{enum:fact_gap.i} and~\eqref{eq:Q(w)_also}~imply
that~$\cP_{L_k}\cup\cP_{L_{k'}}$ is a packing.  Indeed,
let~$w=(w_i)_{1\leq i\leq d}\in L_k$ and
any~$w'=(w_i')_{1\leq i\leq d}\in L_{k'}$ be given.  Then, by
definition, there is some~$i$ such that~$w_i<k<k'=w_i'$.  This implies
that~$Q(w)=Q^{(k)}(w)$ and~$Q(w')=Q^{(k')}(w')$ are disjoint `in the
$i$th dimension' (see Fact~\ref{fact:gap}\eqref{enum:fact_gap.i}).

\begin{definition}[Separated families]
  \label{def:separated_families}
  Let~$\cL=(L_k)_{2\leq k\leq S}$ be a family of $k$-codes
  $L_k\subset[k]^d$.  If, for every~$2\leq k<k'\leq S$, the
  codes~$L_k$ and~$L_{k'}$ are separated, then we say that~$\cL$ is a
  \textit{separated} family of codes.
\end{definition}

\begin{remark}
  \label{rem:warm-up} For~$2\leq k\leq d$,
  let~$L_k=\big\{w=(w_i)_{1\leq i\leq d}\in[k]^d\colon w_k=k\text{ and
  }w_i<k\text{ for all~$i\neq k$}\big\}$.
  Then~$\cL=(L_k)_{2\leq k\leq d}$ is a separated family of gapped
  codes.  Fix~$0<\epsilon\leq d^{-2}$.
  Consider~$\cP=\bigcup_{2\leq k\leq d}\cP_{L_k}$ with~$\cP_{L_k}$
  as in~\eqref{eq:cP_L_k_def}.  Since each~$L_k$ is gapped,
  the~$\cP_k$ are packings.  Also, since~$\cL=(L_k)_{2\leq k\leq d}$
  is a separated family, $\cP$~is a packing.  Furthermore, we
  have~$\nu_k(\cP)=|L_k|=(k-1)^{d-1}$
  {\rm(}recall~\eqref{eq:nu_k(U)_def}{\rm)}
  and~$\val(\cP)=\sum_{2\leq k\leq d}1/(k-1)\sim\log d$
  {\rm(}recall~\eqref{eq:val_def}{\rm)}.  The existence of~$\cP$
  implies a weak form of Theorem~\ref{thm:lwbd_prbsa} {\rm(}namely, a
  lower bound of~$\Omega(\log d)$ instead of~$\Omega(d/\log d)${\rm)}.
\end{remark}

Remark~\ref{rem:warm-up} above illustrates the use we wish to make of
separated families of gapped codes.  Our focus will now shift onto
producing much `better' families than the one explicitly defined in
Remark~\ref{rem:warm-up}.
Indeed, we now prove Lemma~\ref{lem:sep_langs} below, which asserts
the existence of such better families.  We shall need the following
auxiliary lemma.

\begin{lemma}
  \label{lem:cF}
  There is an absolute constant~$d_0$ such that, for any~$d\geq d_0$,
  there are sets~$F_1,\dots,F_d\subset[d]$ such that {\rm(i)}~for
  every~$1\leq k\leq d$, we have~$|F_k|=\lceil d/2\rceil$ and
  {\rm(ii)}~for every~$1\leq k<k'\leq d$, we
  have~$|F_k\cap F_{k'}|<7d/26$.
\end{lemma}
\begin{proof}%
  Let~$r=\lceil d/2\rceil$.  We select each~$F_k$ ($1\leq k\leq d$)
  among the $r$-element subsets of~$[d]$ uniformly at random, with
  each choice independent of all others.  Let~$s=7d/26$.
  Note that, for any~$k\neq k'$, we have~$\EE(|F_k\cap
  F_{k'}|)=r^2/d$.  Let~$\lambda=r^2/d$.  Let
  \begin{equation}
    \label{eq:t_from_s}
    t=s-\lambda\geq s-(d/2+1)^2/d
    \geq{7d\over26}-{1\over d}\({d^2\over4}+d+1\)
    \geq{d\over52}-2\geq{d\over53},
  \end{equation}
  as long as~$d$ is large enough.  We may now apply a Chernoff bound
  for the hypergeometric distribution (see, \textit{e.g.}, \cite{janson00:_random_graph}, 
   Theorem~2.10, inequality~(2.12)) to obtain that 
  \begin{equation}
    \label{eq:hyper_dev}
    \PP(|F_k\cap F_{k'}|\geq s)
    =\PP(|F_k\cap F_{k'}|\geq\lambda+t)
    \leq\exp\(-{2(d/53)^2\over\lceil d/2\rceil}\)
    \leq\ee^{-{3d/53^2}}
  \end{equation}
  for every large enough~$d$.  Therefore, the expected number of
  pairs~$\{k,k'\}$ with $1\leq k<k'\leq d$ for which~$|F_k\cap
  F_{k'}|\geq s$ is less than~$d^2\exp(-3d/53^2)$, 
  which tends to~$0$ as~$d\to\infty$.  Therefore, for any large
  enough~$d$, a family of sets~$F_1,\dots,F_d$ as required does exist.
\end{proof}

We are now ready to state and prove the lemma that asserts the
existence of a separated family of gapped codes that is  `better' than
the one defined in Remark~\ref{rem:warm-up}.

\begin{lemma}[Many large, separated gapped codes]
  \label{lem:sep_langs}
  There is an absolute constant~$d_0\geq2$ such that, for
  any~$d\geq d_0$, there is a separated
  family~$\cL=(L_k)_{2\leq k\leq S}$ of gapped
  $k$-codes~$L_k\subset[k]^d$ such that
   \begin{equation}
    \label{eq:sep_langs.Ls}
    |L_k|\geq{10\over11}(k-1)^d
  \end{equation}
  for every~$2\leq k\leq S$, where
  \begin{equation}
    \label{eq:sep_langs.Q}
    S=\llceil2d\over9\log d\rrceil.
  \end{equation}
\end{lemma}
\begin{proof}
  Let~$S$ be as in~\eqref{eq:sep_langs.Q}
  and let~$F_1,\dots,F_d$ be as in
Fact~\ref{lem:cF}.  In what follows, we only use the~$F_k$ for~$2\leq
k\leq S$.  For each~$2\leq k\leq S$, we 
construct~$L_k\subset[k]^d$ in two parts.  Suppose first that we
have~$L_k'$ with
\begin{equation}
  \label{eq:L_k'.1}
  L_k'\subset([k]\setminus\{k-1\})^{F_k}
  =\{w=(w_i)_{i\in F_k}:w_i\in[k]\setminus\{k-1\}\text{ for all
  }i\in F_k\}.
\end{equation}
We then set
\begin{equation}
  \label{eq:L_k.1}
  \begin{split}
    L_k&=L_k'\times[k-1]^{[d]\setminus F_k}\\
    &=\{w=(w_i)_{1\leq i\leq d}\colon
    \exists w'=(w_i')_{i\in F_k}\in L_k' \text{ such that
    }w_i=w_i'\text{ for all }i\in F_k\\
    &\quad\qquad\qquad\qquad\quad\enspace\,
    \text{ and }w_i\in[k-1]\text{ for all }i\in[d]\setminus F_k\}.
  \end{split}
\end{equation}
Note that, by~\eqref{eq:L_k'.1} and~\eqref{eq:L_k.1}, the
$k$-code~$L_k$ will be gapped ($k-1$ is missed at
every~$i\in F_k$ and~$k$ is missed at every~$i\in[d]\setminus F_k$).
We shall prove that there is a suitable choice for the~$L_k'$
with~$|L_k'|\geq(10/11)(k-1)^{|F_k|}$, ensuring
that~$\cL=(L_k)_{2\leq k\leq S}$ is separated.  Since we shall then
have
\begin{equation}
  \label{eq:L_k.large}
  |L_k|=|L_k'|(k-1)^{d-|F_k|}\geq{10\over11}(k-1)^d,
\end{equation}
condition~\eqref{eq:sep_langs.Ls} will be satisfied and
Lemma~\ref{lem:sep_langs} will be proved.  We now proceed
with the construction of the codes~$L_k'$ ($2\leq k\leq S$).

Fix~$2\leq k\leq S$.  For~$2\leq\ell<k$,
let~$J(\ell,k)=F_k\setminus F_\ell$, and note that
\begin{equation}
  \label{eq:J_large}
  |J(\ell,k)|>\llceil{d\over2}\rrceil-{7\over26}d
  \geq{3\over13}d.
\end{equation}
Let~$v=(v_i)_{i\in F_k}$ be an element
of~$([k]\setminus\{k-1\})^{F_k}$ chosen uniformly at random.  For
every~$2\leq\ell<k$, let us say that~$v$ is \textit{$\ell$-bad}
if~$v_i\neq k$ for every~$i\in J(\ell,k)$.  We have
\begin{equation}
  \label{eq:PP_v_ell-bad}
  \PP(v\text{ is $\ell$-bad})=\(1-{1\over k-1}\)^{|J(\ell,k)|}
  \leq\ee^{-|J(\ell,k)|/S}
  \leq\exp\(-{3d\over13\lceil2d/9\log d\rceil}\)
  \leq d^{-1},
\end{equation}
for every large enough~$d$.  Let us say that~$v$ is \textit{bad} if
it is $\ell$-bad for some~$2\leq\ell<k$.  It follows
from~\eqref{eq:PP_v_ell-bad} that 
\begin{equation}
  \label{eq:PP_v_bad}
  \PP(v\text{ is bad})\leq Sd^{-1}\leq{1\over4\log d}\leq{1\over11}
\end{equation}
if~$d$ is large enough.  Therefore, at least~$(10/11)(k-1)^{|F_k|}$
words~$v\in([k]\setminus\{k-1\})^{F_k}$ are not bad, as long as~$d$
is large enough.  We let~$L_k'\subset([k]\setminus\{k-1\})^{F_k}$ be
the set of such good words.

To complete the proof, it remains to show that the family
$\cL=(L_k)_{2\leq k\leq S}$ is separated. 
More precisely, we show that with the above choice of~$L_k'$
$(2\leq k\leq S)$, the family~$\cL=(L_k)_{2\leq k\leq S}$ with~$L_k$
as defined in~\eqref{eq:L_k.1} is separated.

To this end, fix~$2\leq\ell<k\leq S$.  We show that~$L_\ell$ and~$L_k$
are separated.  Let~$u=(u_i)_{1\leq i\leq d}\in L_\ell$
and~$w=(w_i)_{1\leq i\leq d}\in L_k$ be given.  By the definition
of~$L_k$, there is~$v=(v_i)_{i\in F_k}\in L_k'$ such that~$w_i=v_i$
for all~$i\in F_k$.  Furthermore, since~$v\in L_k'$ is not a bad word,
it is not $\ell$-bad.  Therefore, there
is~$i_0\in J(\ell,k)=F_k\setminus F_\ell$ for which we
have~$v_{i_0}=k$.  Observing that~$i_0\notin F_\ell$ and recalling the
definition of~$L_\ell$, we see that $u_{i_0}<\ell<k=v_{i_0}=w_{i_0}$,
as required.

The  proof of Lemma~\ref{lem:sep_langs} is now complete.
\end{proof}

\subsection{The packing $\,\cU$ in Lemma~\ref{lem:large_value_pack}}
\label{sec:packing-cu_epsilon}
Fix~$\cL=(L_k)_{2\leq k\leq S}$, a separated family of gapped
$k$-codes~$L_k\subset[k]^d$.  We now give, for every sufficiently
small~$\epsilon>0$, the construction of a packing~$\cU_\epsilon(\cL)$
of $d$-hypercubes into the unit bin~$[0,1]^d$ using~$\cL$ and prove
that~$\cU_\epsilon(\cL)$ is indeed a packing.  Choosing~$\cL$ as in
Lemma~\ref{lem:sep_langs} above, we shall deduce
Lemma~\ref{lem:large_value_pack} by taking~$\cU=\cU_\epsilon(\cL)$.

\begin{definition}[Packing~$\cU_\epsilon=\cU_\epsilon(\cL)$]
  \label{def:cU_epsilon(cL)} Suppose $\cL=(L_k)_{2\leq k\leq S}$ is a
  separated family of gapped \hbox{$k$-codes} $L_k\subset[k]^d$.
  Let~$0<\epsilon\leq S^{-2}$.  We put
  \begin{equation}
    \label{eq:5}
    \cU_\epsilon=\cU_\epsilon(\cL)=\bigcup_{2\leq k\leq S}\cP_{L_k},    
  \end{equation}
  where~$\cP_{L_k}$ is as in~\eqref{eq:cP_L_k_def}.
\end{definition}

In Lemma~\ref{lem:indeed_packs} below, we compile the properties that
we need of~$\cU_\epsilon$.  For the relevant notation,
recall~\eqref{eq:K(U)_def}, \eqref{eq:nu_k(U)_def} and
Definition~\ref{def:(1+epsilon)/ZZ}.

\begin{lemma}  %
  \label{lem:indeed_packs}
  Suppose $\cL=(L_k)_{2\leq k\leq S}$ is a separated family of
  non-empty gapped $k$-codes~$L_k\subset[k]^d$.
  Suppose~$0<\epsilon\leq S^{-2}$.
  Let~~$\cU_\epsilon=\cU_\epsilon(\cL)$ be the family of all the
  hypercubes~$Q(w)=Q^{(k)}(w)\subset[0,1]^d$ with~$w\in L_k$
  and~$2\leq k\leq S$.  Then the following assertions hold:
  {\rm(i)}~the hypercubes in~$\cU_\epsilon$ are pairwise disjoint and
  form an $\epsilon$-packing; %
  {\rm(ii)}~for every~$2\leq k\leq S$, we
  have~$\nu_k(\cU_\epsilon)=|L_k|$; %
  {\rm(iii)}~$|K(\cU_\epsilon)|=S-1$. %
\end{lemma}
\begin{proof}
  Let us first check that the hypercubes~$Q(w)$ in~$\cU_\epsilon$ are
  pairwise 
  disjoint.  We remark that, when introducing the notions of gapped
  and separated codes, we already discussed the reason why the~$Q(w)$
  in~$\cU_\epsilon$ are indeed pairwise disjoint.  However, we give a
  formal proof here for completeness.
  Let~$w=(w_i)_{1\leq i\leq d}\in L_k$
  and~$w'=(w_i')_{1\leq i\leq d}\in L_{k'}$
  with~$2\leq k\leq k'\leq S$ be given.  Consider~$Q(w)=Q^{(k)}(w)$
  and~$Q(w')=Q^{(k')}(w')$.  We have to show that
  \begin{equation}
    \label{eq:Qs_disjoint}
    Q(w)\cap Q(w')=\emptyset.
  \end{equation}

  Suppose first that~$k=k'$.  In that case, both~$w$ and~$w'$ are
  in~$L_k=L_{k'}$ and we may suppose that~$w\neq w'$.  Thus, there is 
  some~$1\leq i\leq d$ 
  such that~$w_i\neq w_i'$.  Furthermore, since~$L_k$ is gapped,
  either~$k-1$ or~$k$ is missed by~$L_k$ at~$i$.  In particular, the
  pair~$\{w_i,w_i'\}$ cannot be the pair~$\{k-1,k\}$ and
  therefore
  \begin{equation}
    \label{eq:disjoint_Ii}
    I^{(k)}(w_i)\cap I^{(k)}(w_i')=\emptyset    
  \end{equation}
  (recall Fact~\ref{fact:gap}\eqref{enum:fact_gap.ii}).
  Expression~\eqref{eq:Q(w)_also} applied to~$Q(w)$ and~$Q(w')$,
  together with~\eqref{eq:disjoint_Ii},
  confirms~\eqref{eq:Qs_disjoint} when~$k=k'$.  

  Suppose now that~$k<k'$.  Since~$L_k$ and~$L_{k'}$ are separated,
  there is some~$1\leq i_0\leq d$ such that~$w_{i_0}<k<k'=w_{i_0}'$.
  Fact~\ref{fact:gap}\eqref{enum:fact_gap.i} tells us that
  \begin{equation}
    \label{eq:disjoint_Ii0}
    I^{(k)}(w_{i_0})\cap I^{(k')}(w_{i_0}')=\emptyset.
  \end{equation}
  Expression~\eqref{eq:Q(w)_also} applied to~$Q(w)$ and~$Q(w')$,
  together with~\eqref{eq:disjoint_Ii0},
  confirms~\eqref{eq:Qs_disjoint} in this case also.  We therefore
  conclude that~$\cU_\epsilon$ is indeed a packing.

  The hypercubes in $\cU_\epsilon$ are copies of the
  hypercubes~$Q_k^+$ for~$2\leq k\leq S$, and
  therefore $\cU_\epsilon$ is an $\epsilon$-packing.
  This concludes the proof of Lemma~\ref{lem:indeed_packs}(i).
  Assertions~(ii) and~(iii) are clear.
\end{proof}

We are now ready to prove Lemma~\ref{lem:large_value_pack}.

\begin{proof}[of Lemma~\ref{lem:large_value_pack}]
  Let~$d_0$ be as in Lemma~\ref{lem:sep_langs}.  We may and shall
  suppose that~$d_0\geq e^2$ and that~$d_0$ is large enough so that,
  for every~$d\geq d_0$, the last inequality in~\eqref{eq:4} below
  holds.  We prove that Lemma~\ref{lem:large_value_pack} holds with
  this choice of~$d_0$.  Let~$d\geq d_0$ and~$0<\epsilon\leq d^{-2}$
  be given.  Let $S=\lceil2d/9\log d\rceil$.  Note
  that~$\epsilon\leq d^{-2}\leq S^{-2}$.
  Let~$\cL=(L_k)_{2\leq k\leq S}$ be a separated family of gapped
  $k$-codes as given by Lemma~\ref{lem:sep_langs}.
  Lemma~\ref{lem:indeed_packs} tells us
  that~$\cU_\epsilon=\cU_\epsilon(\cL)$ is an $\epsilon$-packing with
  \begin{multline}
    \label{eq:4}\qquad\quad
    \val(\cU_\epsilon)=\sum_{k\in
      K(\cU_\epsilon)}{\nu_k(\cU_\epsilon)\over(k-1)^d}
    =\sum_{k\in
      K(\cU_\epsilon)}{|L_k|\over(k-1)^d}\\
    \geq{10\over11}(S-1)
    ={10\over11}\(\llceil2d\over 9\log d\rrceil-1\)
    \geq{d\over5\log d}.\qquad\quad
  \end{multline}
  Thus, to prove Lemma~\ref{lem:large_value_pack}, it suffices to 
  take~$\,\cU=\cU_\epsilon$. 
\end{proof}

\section{Concluding remarks}
\label{sec:concluding-remarks}

We have not optimized the numerical constants in our calculations
above.  In particular, the constant~$10$ in
Theorem~\ref{thm:lwbd_prbsa} can be made arbitrarily close to~$4$,
although~$d_0$ would grow as we do so.  We note that, since the
problem posed by \cite{EpsteinS05} is of an asymptotic nature
($d\to\infty$), the specific value of~$d_0$ is not particularly
relevant.

Our approach for finding a certain good packing in this paper is based
on establishing the existence of certain specific families of
compatible codes by the probabilistic method.  We hope similar ideas
will be useful in other related contexts.

\bibliographystyle{abbrvnat}
\bibliography{extracted}
\label{sec:biblio}

\end{document}